\documentclass[11pt,reqno]{amsart}

\newcommand\version{November 16, 2026}

\usepackage{amsmath,amsfonts,amsthm,amssymb,amsxtra}
\usepackage{stmaryrd}
\usepackage{bbm} % to get \1
\usepackage{hyperref}	% to get hyperlinks in equations and references
\usepackage{mathrsfs}
\usepackage{esint}
\usepackage{xfrac}
\usepackage{relsize, pgfplots, tikz, float}
\pgfplotsset{compat=1.18}

\renewcommand{\baselinestretch}{1.1}

\newtheorem{theorem}{Theorem}[section]
\newtheorem{proposition}[theorem]{Proposition}
\newtheorem{lemma}[theorem]{Lemma}

\theoremstyle{definition}

\theoremstyle{remark}

\newtheorem{remark}[theorem]{Remark}

\numberwithin{equation}{section}

\renewcommand{\epsilon}{\varepsilon}

\newcommand{\N}{\mathbb{N}}

\renewcommand{\phi}{\varphi}
\newcommand{\R}{\mathbb{R}}

\newcommand{\Z}{\mathbb{Z}}

\newcommand{\lo}{\mathrm{lo}}

\newcommand{\hi}{\mathrm{hi}}

\usepackage{mathtools, thmtools}

\let\oldtocsection=\tocsection

\let\oldtocsubsection=\tocsubsection

\let\oldtocsubsubsection=\tocsubsubsection

\renewcommand{\tocsection}[2]{\hspace{0em}\oldtocsection{#1}{#2}}
\renewcommand{\tocsubsection}[2]{\hspace{1em}\oldtocsubsection{#1}{#2}}
\renewcommand{\tocsubsubsection}[2]{\hspace{2em}\oldtocsubsubsection{#1}{#2}}

\title[On K\"onig's non-existence conjecture]{On K\"onig's non-existence conjecture \\for the one-dimensional fractional \\Bianchi--Egnell inequality}
\author{Jonas W. Peteranderl}
\address[Jonas W. Peteranderl]
{Mathematisches Institut,
Ludwig-Maximilians-Universit\"at M\"unchen,
Theresienstr. 39, 80333 M\"unchen, Germany}
\email{peterand@math.lmu.de}
\author{Juncheng Wei}
\address[Juncheng Wei]{Department of Mathematics,
The Chinese University of Hong Kong, Shatin, New Territories, Hong Kong}
\email{wei@math.cuhk.edu.hk}
\author{Yuanze Wu}
\address[Yuanze Wu]{School of Mathematics,
Yunnan Normal University, Kunming 650500, China; Yunnan Key Laboratory of Modern Analytical Mathematics and Applications, Kunming, China}
\email{yuanze.wu@ynnu.edu.cn}

\date{\version}

\thanks{\copyright\,2026 by the authors.
This paper may be reproduced, in its entirety, for non-commercial purposes.}

\begin{document}

\begin{abstract}
K\"onig conjectured that the sharp stability constant in the Bianchi--Egnell inequality for the one-dimensional fractional Sobolev inequality equals the local stability constant and that the associated quotient has no minimizer outside the set of optimizers of the fractional Sobolev inequality. We prove both assertions for~$s<1/2$ sufficiently close to~$1/2$, with an explicitly computable threshold. The proof is based on a reinterpretation of the Frank--Lieb rigidity mechanism. Due to the lack of smooth variational structure of the Bianchi--Egnell quotient, we work with a refined limiting inequality of Beckner--Onofri-type on the circle with vanishing first-order moments, which we derive from Toeplitz determinants.
\end{abstract}
\maketitle

\section{Introduction and main results}

For $s\in (0,d/2)$, $d\in \N$, the sharp fractional Sobolev  inequality
  on $\mathbb{R}^d$ is given by
\begin{equation}\label{eq:frac_Sob}
    \|(-\Delta)^{s/2}f\|_{L^2(\R^d)}^2\geq  S_{d,s} \|f\|_{L^{2^*}(\R^d)}^2
\end{equation} with
$2^*=2d/(d-2s)$. The set of optimizers is given by 
\begin{equation*}
  \mathcal{M}_{d,s}\coloneqq
 \{c(a^2+|x-b|^2)^{-(d-2s)/2}:
      c\in\mathbb{R},\ a>0,\ b\in\mathbb{R}^d\} \,. \end{equation*} 
      This follows from the corresponding results of Lieb \cite{Lieb} for the
  Hardy--Littlewood--Sobolev (HLS) inequality by duality (see also Beckner's spherical formulation \cite{Beckner}).

In this paper, we consider the case $d=1$ and write
$S_s=S_{1,s}$ and $\mathcal{M}_s=\mathcal{M}_{1,s}$.
For $f\in\dot H^s(\mathbb{R})\setminus\mathcal{M}_s$, we define the stability quotient and its best constant by
\begin{equation}\label{eq:euclidean}
 \mathcal{E}_s(f)\coloneqq
 \frac{ \|(-\Delta)^{s/2}f\|_{L^2(\R)}^2-S_s \|f\|_{L^{2^*}(\R)}^2}
 {\operatorname{dist}_{\dot H^s}(f,\mathcal{M}_s)^2}\qquad \text{and}\qquad
   c_{\mathrm{BE}}(s)\coloneqq\inf_{f\notin\mathcal{M}_s}\mathcal{E}_s(f)\,.
\end{equation}Bianchi and Egnell~\cite{BE} established stability for the classical Sobolev inequality. The fractional result of Chen, Frank, and Weth~\cite{CFW} implies~$c_{\mathrm{BE}}(s)>0$ in the present range. A second-order computation at one of the optimizers gives the upper bound
\begin{equation}\label{eq:localconstant}
 c_{\mathrm{BE}}(s)\leq c_{\mathrm{loc}}(s)\coloneqq\frac{4s}{3+2s}\,.
\end{equation}
The right side is the asymptotically sharp \emph{local stability constant}, and we refer to \cite[Proposition 2 and its proof]{CFW} for an explicit sequence realizing this constant.

In dimensions at least two, K\"onig proved  that the best constant in the fractional Bianchi--Egnell inequality is strictly smaller than its local counterpart \cite{KoenigSharp} and that  the stability quotient has a minimizer \cite{KoenigExistence}. In dimension one, his local theorem~\cite{KoenigExceptional} instead bounds the quotient $\mathcal{E}_s(f)$ below by $c_{\mathrm{loc }} (s)$ when $f$ is sufficiently close to the set of Sobolev optimizers in the homogeneous Sobolev norm. K\"onig conjectured \cite[Conjecture 2]{KoenigExceptional} that this lower bound holds globally and that the best constant $c_{\mathrm{BE}}(s)$ is never attained outside $\mathcal M_s$ for all $s\in(0, 1/2)$. In this paper, we prove these assertions for $s$ close to the critical value $1/2$ with an explicitly computable lower bound $s_0$.

\begin{theorem}\label{thm}
    There is an explicitly computable $s_0\in(0,1/2)$ such that, for every $s\in (s_0,1/2)$, the quotient  $\mathcal E_s$ has no global minimizer in $\dot H^s(\mathbb R)\setminus\mathcal M_s$ and
$$c_{\mathrm{BE}} (s)= c_{\mathrm{loc}} (s)\,.$$
\end{theorem}

\begin{remark}
The lower bound $s_0$ can be estimated as $s_0\approx 1/2-e^{-10^{55}}$.
\end{remark}

Thus, every admissible function $f$ has a quotient $\mathcal{E}_s(f)$ that is strictly larger than $c_{\mathrm{loc}}$, even though functions $f$ approaching $\mathcal M_s$ have quotients $\mathcal{E}_s(f)$ converging to this value. The quotient $\mathcal{E}_s(f)$ itself is undefined on $\mathcal M_s$. Indeed, the exponent $2^*$ diverges as $s\uparrow 1/2$ and we obtain a Beckner--Onofri-type quotient (see \cite{LebedevMilin,Onofri,Beckner}) as a formal limit.

\subsection{Related results} We emphasize that the existence of optimizers for the stability inequality as in  K\"onig's work \cite{KoenigExistence} is a separate question from  the classification of optimizers for the original inequality. Related works on attainment of stability optimizers concern Caffarelli--Kohn--Nirenberg (CKN) inequalities \cite{WeiWuExistence,DengTianWei}, the Sobolev inequality on the Heisenberg group \cite{TangZhangZhang}, and fractional Sobolev trace inequalities~\cite{ZhangZhouZou}. Attainment has also been studied for the quantitative logarithmic Sobolev inequality \cite{DovettaSerra}, the planar quantitative isoperimetric inequality \cite{BianchiniCroceHenrot}, and the quantitative Yamabe inequality \cite{YamabeExistence}. 

With the goal to understand stability constants better, constructive proofs of quantitative stability have been investigated. In their breakthrough result \cite{DEFFL}, Dolbeault, Esteban, Figalli, Frank, and Loss established explicit positive lower bounds on the stability constants with optimal dimensional dependence for the Sobolev inequality. This has been extended to HLS-inequalities \cite{CLTExplicit} as well as fractional and higher orders \cite{CLTAsymptotic, CLTHigher}. Another approach uses entropy methods giving quantitative Gagliardo--Nirenberg--Sobolev (GNS) and spherical interpolation inequalities with explicit constants \cite{BDNS,BDS}. Recently, explicit sharp stability constants were also obtained for some reverse Sobolev and isoperimetric-type inequalities by direct methods \cite{KoenigReverse,LP25}. 

\subsection{Proof strategy}
Our strategy is to exclude global minimizers of the quotient in \eqref{eq:euclidean} at or below the local level $c_{\mathrm{loc}}(s)$. Compactness below the level $c_{\mathrm{loc}}(s)$ and the local upper bound~\eqref{eq:localconstant} then determine~$c_{\mathrm{BE}}(s)$. We first observe that a potential minimizer must have one sign by the energy comparison in Proposition~\ref{prop:combin} and, after a conformal change and scaling, can be normalized so that~$1$ is a closest Sobolev optimizer. This gives uniform bounds on its bubble correlations and an Euler--Lagrange equation with a probability mixture of closest bubbles; see Section~\ref{sec:2} for details. We emphasize that $d=1$ is the key point in the energy comparison.

Section \ref{sec:3} retains this normalization without assuming balancing. We use the energy bound and the Euler--Lagrange equation to control the $L^\infty$-norm of $u-1$, which yields a uniform bound on $\|\log\rho\|_{H^1}$ for the critical probability density $\rho=u^{1/a}/\int u^{1/a}\,\mathrm d\mu$.  A further conformal change produces a balanced density, with vanishing first moment, preserving this bound.

In Section~\ref{sec:4}, we consider a limiting Beckner--Onofri-type inequality \cite{LebedevMilin,Onofri,Beckner} under the vanishing first moment.  Its energy coefficient drops from~$1/2$ to~$1/4$, with equality densities $Q_z(\theta)=P_z(2\theta)$ (\cite{OPS,Widom} and~\cite[Theorem 6.1, $m=1$]{ChangHang}).  Simon's Toeplitz bound~\cite[Theorems 1.1 and 2.4]{Simon} then yields an entropy remainder measuring distance to the~$Q_z$ with matching second moment. To quantify this remainder, we expand the auxiliary functional $\Phi_a$ defined by \eqref{eq:Phi} and use $\Phi_a(\rho^a)\leq 0$, which gives $\|\rho/Q_z-1\|_{\infty}=O(a^{1/4})$.  Linearizing the remainder further yields coercivity for constrained odd perturbations.

A low-high frequency decomposition transfers the odd coercivity established in Section~\ref{sec:4} to positive $a$, which, together with balancing, controls the moment error, while parity cancels the linear norm variation.  To further control the even perturbations, we introduce the technique of angle doubling in Section~\ref{sec:5}, which proves nonnegativity on the even part of $\rho^a$.  This leads to $\Phi_a(\rho^a)>0$ unless $\rho=1$. The key point in this section is the use of a constrained second-variation argument related to~\cite{FrankLieb,CaseFL}, which simplifies the argument by avoiding further differentiation of the supremum over bubble correlations.

In Section \ref{sec:6}, we choose the threshold. Indeed, $\Phi_a(\rho^a)>0$ contradicts the quotient bound $\Phi_a(\rho^a)\leq0$, while $\rho=1$ would make the original minimizer a Sobolev optimizer. Compactness below the local level and a local test sequence then identify the best stability constant.

\subsection*{Acknowledgments}
J. W. Peteranderl acknowledges support from the Deutsche For\-schungsgemeinschaft (DFG; TRR~352, Project-ID~470903074, and FR~2664/3-1), the Studienstiftung des deutschen Volkes, and the Fondation Sciences Math\'ematiques de Paris (FSMP). J. Wei was supported by National R\&D Program of China (Grant No.~2022YFA1005602), NSFC (Nos.~12631008), and Hong Kong General Research Fund (No.~14303125) ``On Fujita equation in the critical or supercritical regime''. Y. Wu acknowledges support from NSFC (Nos.~12631008 and 12671144), the Yunnan Provincial Innovation Team on the Interdisciplinary Integration of Modern Applied Mathematics and Life Sciences (No.~202405AS350003), the Yunnan Fundamental Research Projects (No.~202601CJ070001), and the Yunnan Revitalization Talent Support Program.

The authors thank 
the Tianyuan Mathematics Research Center in Kunming for its hospitality at the workshop ``Functional inequalities, their stability and related PDEs'' in October 2026. J. W. Peteranderl also thanks Yunnan Normal University and Yuanze Wu for their hospitality through the scientific program week in November 2026.

\subsection*{AI declaration}
This work was developed with assistance from GPT~6 Astra. The authors supplied initial notes for a proof based on an exclusion strategy using the Frank--Lieb mechanism \cite{FrankLieb}. They guided the investigation with particular emphasis on proving nonattainment without symmetry assumptions on the competing functions and on obtaining a computable threshold $s_0$.

An initial proof strategy used commutator and second-variation estimates applied directly to the Bianchi--Egnell quotient. This led to technically involved calculations because the distance term contains a supremum over bubble correlations, 
which may be nonsmooth when several bubbles are equally close; cf.~\eqref{eq:circleQ}. The available control was already perturbative at this stage.

The authors simplified this approach by introducing a quantitative remainder for the constrained Beckner--Onofri inequality arising in the limit $s\uparrow 1/2$. This avoids differentiating the supremum over bubble correlations and reduces the problem, after uniform estimates and conformal balancing, to a perturbative comparison with the equality family of the limiting inequality. A further simplification introduced by the authors is the argument in Section~\ref{sec:5}: an antipodal even--odd decomposition makes the relevant linear variation vanish by parity, while the odd component is controlled by a constrained quadratic spectral estimate.

The AI assisted in developing calculations and draft arguments and with literature searches and mathematical checks. The authors remain responsible for the mathematical claims and their presentation in writing.

\section{Preliminaries}\label{sec:2}

We write $\lesssim_R$ if we have $\leq$ up to a constant depending on some quantity $R$; we use $\gtrsim_R$ analogously. In particular, constants may change from line to line.  Throughout the following sections, we set
\begin{equation} \label{eq:begin}
 a \coloneqq \frac12-s \in \left(0, \frac12\right)
 \qquad \text{and} \qquad c_{\mathrm{loc}} \coloneqq  c_{\mathrm{loc}}(s)= \frac{ 1-2a}{2-a}\, .
\end{equation}

\subsection{Passing to the circle}
In the following, we work on the circle with normalized arclength.  To this end, we parametrize $\mathbb{S}^1\subset\mathbb{R}^2$ by $\omega(\theta)=(\cos\theta,\sin\theta)$, and set
\begin{equation*} \mathrm d\mu (\theta) \coloneqq\frac{\mathrm d\theta}{2\pi}\qquad \text{and}\qquad \mathbb D\coloneqq\{z\in{\mathbb R}^2:|z|<1 \}\,
.
\end{equation*}
Define the normalized conformal fractional Laplacian~$T_a$ on $H^{1/2-a}(\mathbb S^1)$ by $T_ae^{in\theta}=\lambda_{|n|}(a)e^{in\theta}$, where 
\begin{equation}\label{eq:symbol}
 \lambda_0(a)\coloneqq1\qquad\text{and}\qquad
 \lambda_n(a)\coloneqq
 \frac{\Gamma(n+1-a)\Gamma(a)}{\Gamma(n+a)\Gamma(1-a)}
\end{equation} for $n\geq 1$. We write $\lambda_n=\lambda_n(a)$ if the parameter is clear from the context. 

We denote the norm and inner product with respect to~$\mu$ by $\|\cdot\|_p$ and $\langle\cdot,\cdot\rangle$, respectively. The sharp fractional Sobolev inequality \eqref{eq:frac_Sob} in one dimension becomes \begin{equation}\label{eq:ineq}
    \langle u,T_au\rangle\geq\|u\|_{1/a}^2\,.
\end{equation}
Next, we transform the optimizers. We
 define the \emph{Poisson kernel} by
\begin{equation}\label{eq:bubbles}
 P_z(\theta) \coloneqq\frac{1-|z|^2}{1- 2 z\cdot\omega(\theta)+|z|^2}\,,\qquad z\in\mathbb{D} \,,   \end{equation} which satisfies $T_a(P_z^a)=P_z^{1-a}$ and $
 \langle P_z^a,T_a(P_z^a)\rangle=\|P_z^a\|_{1/a}^2=1$. For $0\leq r<1$, we also write $P_r=P_{(r,0)}$.
 The set of optimizers, including zero, is then given by
\begin{equation*}
\mathcal{M}\coloneqq\{tP_z^a :t\in\mathbb{R}, \ z\in\mathbb{D}\} \,.
\end{equation*} In the following, we also use the  norm conventions \begin{equation*}  \|v\|_{H^t}^2 \coloneqq \sum_n(1+n^2)^t| \widehat v(n)|^2  \qquad\text{and}\qquad \|v\|_{C^k}\coloneqq\max_{0\leq j\leq k}\|v^{(j)}\|_\infty\,,
\end{equation*}
where
\begin{equation*}
 \widehat v(n)\coloneqq\int_{\mathbb{S}^1}v(\theta)e^{-in\theta}\,\mathrm d\mu\,.
\end{equation*}

For the regularity and $L^\infty$-estimates in Section \ref{sec:3}, we also need to treat the  inverse of $T_a$. It can be written as a convolution $$T_a^{-1}f(\theta)=\int_{\mathbb S^1}k_a(\theta-\varphi)f(\varphi)\,\mathrm d\mu(\varphi)$$ with the kernel \begin{equation}\label{eq:inversekernel}
 k_a(\theta)\coloneqq \frac1{D_a| 2\sin(\theta/2)|^{2a}} \qquad \text{with}\quad D_a\coloneqq\frac{\Gamma(1-2a)} {\Gamma(1-a)^2}\,.
\end{equation} Indeed, this is a consequence of the Funk--Hecke formula
\begin{equation}\label{eq:funkhecke}
 \int_{\mathbb S^1}\frac{Y_n(\omega')}{|\omega-\omega'|^{2a}}\,\mathrm d\mu(\omega')
 =\frac{D_a}{\lambda_n(a)}Y_n(\omega)
\end{equation} for a spherical harmonic $Y_n$ of degree $n$ on $\mathbb S^1$;  see \cite[Lemma 6]{FrankLectures} with $d=1$, $\alpha=a$, and normalized arclength after applying the duplication formula for Gamma functions. Note that spherical harmonics of degree zero are the constants and of degree $n\geq1$ are spanned by $\cos(n\theta)$ and $\sin(n\theta)$.

Stereographic projection transforms the Euclidean quotient~\eqref{eq:euclidean} into
\begin{equation}\label{eq:circleQ}
 \mathcal{E}_a(u)\coloneqq
 \frac{\langle u, T_au\rangle-\|u\|_{1/a}^2}{\langle u, T_au\rangle-\sup_{z\in\mathbb{D}}|\langle u,P_z^{1-a}\rangle|^2}\,;
\end{equation} cf.~\cite{FrankLectures}. The denominator is the squared distance to the set of Sobolev optimizers in the Hilbert norm $v\mapsto \langle v,T_av\rangle^{1/2}$. Indeed, for each fixed center~$z$, minimizing over the amplitude gives
\begin{align*}
\inf_{t\in\mathbb{R}}\langle u-tP_z^a,T_a(u-tP_z^a)\rangle &=\inf_{t\in\mathbb{R}} (\langle u,T_au\rangle   -2t\langle u,P_z^{1-a}\rangle+t^2 )\\
 &=\langle u,T_au\rangle-|\langle u,P_z^{1-a}\rangle|^2\,.
\end{align*}
Taking the infimum over $z$ proves the asserted formula. By a \emph{conformal transformation of the circle} we mean the boundary restriction~$\varphi$ of a M\"obius automorphism of the unit disk; see \cite[Appendix A]{FPR} or \cite{Ratcliffe}, for instance. Let $J_\varphi$ denote the Jacobian with respect to $\mu$. Then conformal changes of the form
$u\mapsto J_\varphi^a(u\circ\varphi)$ preserve both sides of the inequality \eqref{eq:ineq} and the quantities in the quotient
\eqref{eq:circleQ}. For a nonzero, nonnegative~$u$, we define  its \emph{critical density} as
\begin{equation}\label{eq:criticaldensity}
 \rho\coloneqq \frac{u^{1/a}}{\int_{\mathbb S^1}u^{1/a}\,\mathrm d\mu}\,.
\end{equation} As a conformal change sends  this density $\rho$ to $J_\varphi (\rho\circ\varphi)$, the unit mass of $\rho$ is also preserved under conformal change. We call a probability density $\rho$ \textit{balanced} if $\int_{\mathbb S^1}\rho\omega(\theta)\,\mathrm d\mu(\theta)=0$. Proposition \ref{prop:regularity} obtains the balancing condition by a further conformal change after proving uniform regularity.

Every nonzero $u\in H^{1/2-a}$ has a closest bubble by a standard coercivity argument; see
\cite[Proof of Proposition 2]{CFW}. Indeed, some bubble correlation
is nonzero, so the distance to $\mathcal M$ is strictly less than $\langle u, T_a u\rangle$. Moreover, the bubble correlations tend
uniformly to zero as $|z|\to1$, by approximation with bounded
functions and H\"older's inequality. Thus, all maximizing
parameters lie in a compact subset of $\mathbb D$.

\subsection{Fixed-parameter regularity} 
The $L^\infty$-estimate in Section \ref{sec:3}  uses the integral representation of
the equation at a maximum point of $u$. The following lemma
provides a continuous representative, ensuring that such a
point exists. Recall that $T_a^{-1}$ is convolution with the kernel $k_a$ in~\eqref{eq:inversekernel}. Taking $n=0$ in~\eqref{eq:funkhecke} gives $\|k_a\|_1=1$.

\begin{lemma}\label{lem:weakregularity}
Let $0<a<1/2$ and let $u\in H^{1/2-a}$ solve
\begin{equation}
    \label{eq:equation_nonlinear}cT_au=\lambda|u|^{1/a-2}u-f
\end{equation} weakly, with $c>0$, $\lambda\in\mathbb{R}$, and $f$ smooth.
Then $u$ has a continuous representative on $\mathbb S^1$.
\end{lemma}
\begin{proof}
Sobolev embedding gives $u\in L^{1/a}$, and hence
\begin{equation*}
 V\coloneqq \lambda c^{-1}|u|^{1/a-2}
 \in L^{1/(1-2a)} \,.
\end{equation*} We improve integrability of $u$ using a regularity-lifting theorem as can be found in
\cite[Theorem 6.1]{ChenGuoLi}. Split $$V=V_1+V_2\,,$$ with $V_1$
bounded and $\|V_2\|_{1/(1-2a)}$ sufficiently small. Then 
$$u=T_a^{-1}(V_1u-c^{-1}f)+T_a^{-1}( V_2u)\,.$$
Suppose $u\in L^{r_0}$,
where $r_0\geq1/a$.  Choose finite $r>r_0$ with
$1/r_0-1/r<1-2a$. Let $\ell\in \{r,r_0\}$. Then the HLS-inequality gives
\begin{align*}
 \|T_a^{-1}(V_1u-c^{-1}f)\|_r
& \lesssim_{a,r,r_0} 
\|V_1\|_\infty\|u\|_{r_0}
       +c^{-1}\|f\|_{r_0}<\infty\,,\\\|T_a^{-1}( V_2h) \|_\ell
& \lesssim_{a,\ell} \|V_2\|_{1/(1-2a)}\|h\|_\ell\,, \qquad h\in L^\ell\,.
\end{align*}
Choosing 
$\|V_2\|_{1/(1-2a)}$ small enough, the linear map $h\mapsto T_a^{-1}(V_2h)$ is a
contraction on both $L^{r_0}$ and $L^r$. The lifting theorem yields $u\in L^r$, so iteration reaches 
every finite exponent.

 In particular, we have
\begin{equation*}
 F\coloneqq c^{-1} (\lambda|u|^{1/a-2}u-f ) \in L^{2/(1-2a)},\qquad u=k_a*F\,.
\end{equation*} Since $k_a\in L^{2/(1+2a)}$, H\"older's inequality
makes the convolution bounded and gives
\begin{equation*}
 \|(k_a*F)(\cdot+h)-k_a*F\|_\infty \leq \|k_a(\cdot+h)-k_a\|_{2/(1+2a)}
 \|F\|_{2/(1-2a)} \rightarrow 0
\end{equation*} as $h\to 0$. Thus, $k_a*F$ is a  continuous representative of $u$.
\end{proof}

\begin{remark}
  
    In fact,  we even have $u\in H^{2-2a}(\mathbb S^1)$. To see this, use $T_a^{-1}: H^t \to H^{t+1-2a}$,  which follows from  the asymptotic behavior of the eigenvalues \eqref{eq:symbol}. Since $u$ is bounded, the  nonlinearity  $t\mapsto|t|^{1/a-2}t$ is Lipschitz continuous on bounded sets. The fractional difference  characterization of $H^t$, $0<t<1$, (see \cite[Proposition 3.4]{DNPV12}, for instance) therefore permits bootstrapping the equation \eqref{eq:equation_nonlinear} until $u\in H^1$. The chain rule for $H^1$-functions and one further application of $T_a^{-1}$ yield $u\in H^{2-2a}$. However, this stronger regularity is not needed below.
\end{remark}

\subsection{Sign and the Euler--Lagrange equation}
We will use the following integral kernel representation to exclude sign changes in the next proposition. For smooth $u$, the integral representation of the beta function gives the formula
\begin{equation}\label{eq:Tkernel}
 (T_au)(\theta)
 =u(\theta)+\tau_a\int_{\mathbb{S}^1}
 K_a(\theta-\varphi) (u(\theta)-u(\varphi) )
 \,\mathrm d\mu(\varphi)\,,
\end{equation}
where
\begin{equation*}
 K_a(t)\coloneqq|2\sin(t/2)|^{-2+2a}\qquad \text{and}\qquad \tau_a\coloneqq\frac{2(1-2a)}{a\displaystyle\int_{\mathbb{S}^1}|2\sin(\theta/2)|^{2a}\,\mathrm d\mu}\,.
\end{equation*} The integral is absolutely convergent for $a>0$ as its integrand is of order $O(|\theta-\varphi|^{-1+2a})$ for $|\theta-\varphi|\ll1$. 

The distance to the set of Sobolev optimizers does not have to be differentiable since several bubbles may be equally close. Hence, the closest bubble need not be unique, and the Euler--Lagrange equation involves a probability mixture
of closest bubbles -- as we will see in the next proposition. We denote the set of closest-bubble parameters by 
\begin{equation*} Z(u)\coloneqq \operatorname*{arg\,max}_{z\in\mathbb D}
 |\langle u,P_z^{1-a}\rangle|^2\,,
\end{equation*} so $z\in Z(u)$ if some scalar multiple of $P_z^a$ minimizes the distance to $u$.

\begin{proposition}\label{prop:combin}
A global minimizer $u$ of the quotient in~\eqref{eq:circleQ}, with $u\in H^{1/2-a}\setminus\mathcal{M}$ and
$$\gamma\coloneqq \mathcal{E}_a(u)\leq c_{\mathrm{loc}}\,,$$ has one sign. After a sign change, a conformal transformation, and a multiplication by a positive constant, if necessary, $u$ satisfies
\begin{equation}\label{eq:cap} u\geq 0\,,
 \qquad \| u\|_1=1\,,
 \qquad\text{and}\qquad \langle u,P_z^{1-a}\rangle\leq1\quad \text{for}\ z\in\mathbb{D}\,.
\end{equation}
There is a probability measure $\tau$ supported on the set $Z(u)=\{z\in\mathbb D:\langle u,P_z^{1-a}\rangle=1\}$ such that
\begin{equation}\label{eq:combin}
 (1-\gamma)T_au+\gamma T_ag=\|u\|_{1/a}^{2-1/a} u^{1/a-1}\,,
\end{equation} where \begin{equation}
    \label{eq:g}g\coloneqq \int_{\mathbb D} P_z^a\,\mathrm d\tau(z)
\end{equation} is the probability mixture of closest bubbles.
\end{proposition}

\begin{proof}
We divide the proof into several steps for the sake of clarity.

\emph{Step 1: Excluding sign changes.}
Symmetrizing the integral representation~\eqref{eq:Tkernel}
gives
\begin{equation}\label{eq:quadraticform}
 \langle v,T_av\rangle
 =\|v\|_2^2+\frac{\tau_a}{2}
 \int_{\mathbb S^1}\int_{\mathbb S^1}
 K_a(\theta-\varphi)
  (v(\theta)-v(\varphi) )^2
 \, \mathrm  d\mu(\theta)\,\mathrm  d\mu(\varphi)\,.
\end{equation} This formula holds a priori for smooth $v$ and extends 
to $v\in H^{1/2-a}$ by approximation.
Note that also 
$|v|\in H^{1/2-a}$.

Define
\begin{align*}
 B(v)&\coloneqq
 \sup_{z\in\mathbb D}|\langle v,P_z^{1-a}\rangle|^2\,,\\
 F_\gamma(v)&\coloneqq
  (1-\gamma)\langle v,T_av\rangle- \|v\|_{1/a}^2 +\gamma B(v)\,.
 \end{align*} For $v\not \in \mathcal M$, we have $F_\gamma(v)=(\langle v,T_av\rangle- B(v)) (\mathcal E_a(v)-\gamma)$. 
Global minimality of $u$ gives $F_\gamma(v)\geq0$ for $v\notin\mathcal M$, while $F_\gamma$ vanishes on
$\mathcal M$ and at~$u$. By conformal invariance, we may assume that $z=0$ maximizes $B(|u|)$.  If we write $u=u_+-u_-$, $u_+,u_-\geq 0$, then subtracting \eqref{eq:quadraticform} for $|u|$ from that
for $u$ gives
\begin{align}\notag
 \langle u,T_au\rangle-\langle |u|,T_a|u|\rangle  &=4\tau_a\int_{\mathbb S^1}\int_{\mathbb S^1} K_a(\theta-\varphi)u_+(\theta)u_-(\varphi)
 \,\mathrm d\mu(\theta)\,\mathrm   d\mu(\varphi)
 \\
   &\geq\frac{2(1-2a)}a \|u_+\|_1\|u_-\|_1\,, \label{eq:abs1}
\end{align} 
where we used $\tau_a\min K_a\geq(1-2a)/(2a)$.  As $z=0$ maximizes $B(|u|)$ and is a competitor for $B(u)$, we obtain 
\begin{equation}
    \label{eq:abs2}B(|u|)-B(u)\leq (\|u_+\|_1+\|u_-\|_1)^2-(\|u_+\|_1-\|u_-\|_1)^2= 4\|u_+\|_1\|u_-\|_1\,.
\end{equation}
    Combining \eqref{eq:abs1} and \eqref{eq:abs2} yields
\begin{equation*}
 F_\gamma(|u|)-F_\gamma(u)
  \leq\frac{2\|u_+\|_1\|u_-\|_1}{a} (\gamma-(1-2a) )\leq 0
\end{equation*}
as $\gamma\leq c_{\mathrm{loc}}<1-2a$. If $u$ changed sign, then $\|u_+\|_1,\|u_-\|_1>0$ would contradict $F_\gamma(|u|)\geq0=F_\gamma(u)$.
Thus $u$ has one sign.

Replacing $u$ by $-u$ if necessary, we may assume $u\geq0$. Normalize a closest bubble to the constant bubble and scale $u$ so that $\int u\,\mathrm d\mu=1$. Then $B(u)=1$ and~\eqref{eq:cap} holds.

\emph{Step 2: Deriving the Euler--Lagrange equation.}
Note that the set of the closest-bubble parameters
 $Z(u)=
 \{z\in\mathbb D:\langle u,P_z^{1-a}\rangle=1\}$ is compact. For each fixed direction $v$, all maximizing parameters for $u+tv$, with $|t|$ sufficiently small, lie in one compact subset of $\mathbb D$. Indeed, the correlations of $u$ vanish uniformly at the boundary, and $$|\langle t v,P_z^{1-a}\rangle|\leq |t| \langle v,T_a v\rangle^{1/2}$$ uniformly in $z$. Danskin's theorem~\cite{Danskin}, applied to $Z(u)\subset \mathbb D$, therefore gives the one-sided derivative 
$$ D^+_v B(u)\coloneqq  \lim_{t\downarrow 0} \frac{B(u+tv)-B(u)}t=2\max_{z\in Z(u)}
\langle P_z^a ,T_av \rangle\,. $$ Since $F_\gamma$ attains its minimum at $u$, its one-sided
directional derivative satisfies \begin{equation*}
 \frac 12 D^+_vF_\gamma (u)= \left\langle (1-\gamma)u-\|u\|_{1/a}^{2-1/a}T_a^{-1}(u^{1/a-1}),T_av
 \right\rangle +\gamma \max_{z\in Z(u)} \langle P_z^a,T_av\rangle \geq0\,.
\end{equation*}
 Since $\gamma=\mathcal E_a(u)>0$ by stability of the fractional Sobolev inequality, the Hahn--Banach theorem implies that \begin{equation}\label{eq:elem}
 \gamma^{-1}(\|u\|_{1/a}^{2-1/a}T_a^{-1}(u^{1/a-1})-(1-\gamma)u) 
\in \overline{\operatorname{conv} } \{P_z^a:z\in Z(u)\}\,.
\end{equation}
Otherwise, a separating direction would contradict the
preceding inequality.

As $z \mapsto P_z^a$ is continuous in   $H^{1/2-a}$ on the compact set $Z(u)$,  every element of this closed convex hull can be written as the  integral of these bubbles with respect to a probability measure $\tau$ on $Z(u)$. Denoting the representing integral in \eqref{eq:elem} by $g$ and applying $T_a$ leads to \eqref{eq:combin}. \end{proof}

\section{Uniform bounds for a potential minimizer}\label{sec:3}
For the estimates concerning $u$ in this section, we retain the normalization \eqref{eq:cap} from Proposition \ref{prop:combin}, that is, $u\geq 0$, $\|u\|_1=1$, and all bubble correlations are at most one. We also assume $\mathcal E_a(u)\leq c_{\mathrm{loc}}$, with $c_{\mathrm{loc}}$ as in \eqref{eq:begin}. Lemma \ref{lem:peak} and Proposition \ref{prop:regularity} additionally require \eqref{eq:combin}, with its probability mixture of closest bubbles. Balancing is not assumed but obtained at the end of Proposition \ref{prop:regularity}.
The uniform estimates are obtained in three stages. Bubble correlations first bound $\langle u,T_a u\rangle$. The Euler--Lagrange equation then prevents concentration of $\|u\|_{1/a}$. Finally, the inverse operator yields regularity of the critical density, uniformly in $a$. In some lemmas below, we only need that $a$ is uniformly bounded away from $1/2$; for convenience we impose the uniform upper bound $a\le1/4$.

\subsection{Energy from bubble correlations}
A correlation bound against every bubble also bounds every Poisson average of~$u$.  We compare~$u$ with these averages at the scale $t=3e L^{-1/a}$ associated with the superlevel set $\{u>L\}$. 
This gives an upper bound for the critical integral $\int u^{1/a}\,\mathrm d\mu$ in terms of the energy $\langle u,T_au\rangle$, while $\mathcal E_a(u)\leq c_{\mathrm{loc}}$ gives a lower bound for the same integral. Their incompatibility at large energy proves the following estimate; see \eqref{eq:quadrlin}.

\begin{lemma}\label{lem:cap}
Let $0<a\leq1/4$ and let $u\in H^{1/2-a}\setminus\mathcal M$ satisfy \eqref{eq:cap} and $\mathcal E_a(u)\leq c_{\mathrm{loc}}$. Then we have
\begin{equation}\label{eq:Ta _bound}
 \langle u,T_au\rangle-1\lesssim a\qquad \,.
\end{equation}

\end{lemma}

\begin{proof}
We divide the proof into several steps for the sake of clarity.

\emph{Step 1: Eigenvalue bounds.} We first show that the  eigenvalues $\lambda_n=\lambda_n(a)$ of $T_a$ defined in~\eqref{eq:symbol} obey
\begin{equation} \label{eq:symbolbounds}  e^{-3a}n^{1-2a} \leq a\lambda_n(a)\leq n\,. \end{equation}
for $n\geq 1$. Indeed, we have 
\begin{equation} \label{eq:gammaprod} a\lambda_n=(1-a)\prod_{j=1}^{n-1} \frac{j+1-a}{j+a}\,, 
\end{equation} and the upper bound follows by comparing the factors  in \eqref{eq:gammaprod} with their values at zero. 

To prove the lower bound, set $R_n(a)\coloneqq \log(a\lambda_n/n^{1-2a})$. We compute
\begin{align*}
R_n'(a) &= 2\log n-\sum_{j=1}^{n}\frac1{j-a} -\sum_{j=1}^{n-1}\frac1{j+a}\,,\\
    R_n''(a) &= -\frac1{(1-a)^2}  +\sum_{j=1}^{n-1}\left(\frac1{(j+a)^2}
  -\frac1{(j+1-a)^2}\right)\,.
\end{align*} Evaluation at $0$ gives
$R_n(0)=0$ and $R_n'(0)\geq -2$. The latter bound follows from the elementary bound $1+\log n\geq \sum_{j=1}^n j^{-1}$.  As  $0\leq a<1/2$, we also find
$R_n''(a)\geq - (1-a)^{-2}$.
By Taylor's theorem, we then obtain
\begin{align*}
R_n(a) &
=aR_n'(0) +\int_0^a (a-t)R_n''(t) \,\mathrm dt\geq-2a -\int_0^a\frac{a-t}{(1-t)^2}\,\mathrm dt
\\&\geq-2a-\frac1{(1-a)^2} \int_0^a(a-t)\,\mathrm dt
=-2a-\frac{a^2}{2(1-a)^2}\geq -3a\,,\end{align*} which proves the lower bound in \eqref{eq:symbolbounds}.

\emph{Step 2: Distribution functions.}
The Poisson semigroup $P_t$ acts on functions $u$ via $$(P_tu)(\theta)
=\int_{\mathbb S^1}P_{e^{-t}}(\theta-\varphi)u(\varphi)\,\mathrm d\mu(\varphi)\,.$$ For $0<t\leq1$, the correlation bound in \eqref{eq:cap} and $\|P_{e^{-t}}\|_\infty= \coth (t/2)\leq 3/t$ give $$ \| P_tu\|_\infty \leq \|P_{e^{-t}}\|_\infty^a \int_{\mathbb S^1}P_{e^{-t}}^{1-a}(\theta-\varphi)u(\varphi)\,\mathrm d\mu(\varphi)\leq\left (\frac 3t\right)^a\,. $$ Choose $t\coloneqq 3eL^{-1/a}$ with $L\geq(3e)^a$. Then $P_tu\leq e^{-a}L$, and thus we obtain
\begin{equation}  \mu(\{u>L\}) \leq \mu(\{ |u-P_tu| >(1-e^{-a})L\})
\leq\frac{\|u-P_tu\|_2^2}
 {(1-e^{-a})^2L^2}\label{eq:Cheb}\end{equation} by Chebyshev's inequality.

In a layer-cake formulation of $\|u\|_{1/a}$, we split the integral on $\R_+$ at $(3e)^a$ and apply $\mu(\{u>L\})\leq 1$ and \eqref{eq:Cheb}, respectively, to obtain
\begin{align}\notag
 \int_{\mathbb S^1}u^{1/a}\,\mathrm d\mu &=\frac1a\int_0^\infty L^{1/a-1}\mu(\{u>L\})\,\mathrm dL\\
 &\leq3e+\frac1{a(1-e^{-a})^2}
  \int_{(3e)^a}^\infty L^{1/a-3}\|u-P_tu\|_2^2\,\mathrm dL\,. \label{eq:layer}\end{align} By Parseval's identity,
Tonelli's theorem, and substitution with $y=t|n|=3e|n|L^{-1/a}$, the latter integral equals
\begin{equation}\label{eq:Pars}
    a\sum_{n\ne0}|\widehat u(n)|^2
(3e|n|)^{1-2a}
\int_0^{|n|}y^{-2+2a}(1-e^{-y})^2\,\mathrm dy\,,
\end{equation} The integral in $y$ is uniformly bounded in $n\in \N$ and $0<a\leq 1/4$: use  $1-e^{-y}\leq y$ on $(0,1)$ and $1-e^{-y}\leq1$ on $(1,\infty)$.

\emph{Step 3: Closing the energy bound.}
By~\eqref{eq:symbolbounds}, \eqref{eq:layer}, \eqref{eq:Pars} and
$1-e^{-a}\gtrsim a$, there are computable $C_0,C_1>0$ such that $x\coloneqq a^{-1}(\langle u,T_a u\rangle-1)$ satisfies
\begin{equation}
    \|u\|_{1/a}^{1/a}\leq C_0+C_1 x\,.
\label{eq:upper_crit}
\end{equation}
On the other hand, $\mathcal E_a(u)\leq c_{\mathrm{loc}}$ yields
$$ \|u\|_{1/a}^2\geq 1+\frac{1+a}{2-a}\,ax \geq1+\frac{ax}{2}\,. $$ As $a\mapsto(1+ax/2)^{1/(2a)}$ is non-increasing for fixed $x\geq 0$, we bound
\begin{equation}\label{eq:quadrlin}
    \left(1+\frac18 x\right )^2 \leq \left(1+\frac 12 ax\right)^{1/(2a)}
\leq \|u\|_{1/a}^{1/a}\leq C_0+C_1x\,, 
\end{equation} which provides a computable upper bound on $x$.
\end{proof}

\subsection{Controlling the $L^\infty$-norm}
The proof of Lemma~\ref{lem:cap} also gives $\int_{\mathbb S^1}u^{1/a}\,\mathrm d\mu\lesssim1$, by combining \eqref{eq:Ta _bound} and \eqref{eq:upper_crit}. We now use equation \eqref{eq:combin} to 
strengthen this integral bound to the uniform estimate.

For a function $f$ on $\mathbb S^1$, we denote the discrete Laplacian by
$$\Delta_hf(\theta)\coloneqq f(\theta+h)+f(\theta-h)-2f(\theta)\,.$$

\begin{lemma}\label{lem:kernelbounds}
For $0<a\leq1/4$ and $ |h|\leq 1/10$, the kernel $k_a$ and its normalization $D_a$, defined in \eqref{eq:inversekernel}, satisfy
\begin{equation}
 0\leq\log D_a\lesssim a^2\qquad\text{and}\qquad
 \|\Delta_hk_a\|_1\lesssim a|h|^{1-2a}
 \,.\label{eq:logD_aDelh}
\end{equation}
\end{lemma}

\begin{proof}
Set $F(t)\coloneqq \log\Gamma(1-t)$. Since $F''$ is positive
and bounded on $[0,1/2]$, we can bound
$$
0\leq\log D_a
=\int_0^a\int_0^a F''(s+t)\,\mathrm ds\,\mathrm dt
\lesssim a^2\,.
$$
Thus $D_a\geq1$. Differentiating the kernel $k_a$ on $0<|t|\leq\pi$ gives
$|k_a'(t)|\lesssim a|t|^{-1-2a}$ and $|k_a''(t)|\lesssim a|t|^{-2-2a}$. 

By symmetry, it suffices to consider $0<h\leq1/10$. For $|\theta|<2h$, subtracting the constant $k_a(h)$ from each term of $\Delta_hk_a$ and translating gives
$$
\int_{|\theta|\leq2h}|\Delta_hk_a(\theta)|\,\mathrm d\mu(\theta)
\lesssim\int_0^{3h}|k_a(t)-k_a(h)|\,\mathrm dt\lesssim ah^{1-2a}\,,
$$ where we used the fundamental theorem of calculus and the bound on $|k_a'|$ in the last step. For $|\theta|>2h$, every point between $\theta-h$ and
$\theta+h$ has angular distance at least $|\theta|/2$
from zero. Taylor's formula therefore gives an estimate of $|\Delta_hk_a|$ in terms of $|k''_a|$, and thus
$$
\int_{|\theta|>2h}|\Delta_hk_a(\theta)|\,\mathrm d\mu(\theta)
\lesssim ah^2\int_{2h}^{\pi}t^{-2-2a}\,\mathrm dt
\lesssim ah^{1-2a}\,,
$$ which completes the proof.
\end{proof}

With these preparations, we can now deduce the uniform bound on $u$.

\begin{lemma}\label{lem:peak} There is a computable 
$a_1\in(0,1/4)$ such that, for $0<a\leq a_1$, every $u$ satisfying the assumptions of Lemma \ref{lem:cap} and the Euler--Lagrange equation \eqref{eq:combin}, with $\gamma=\mathcal E_a(u)$ and the probability mixture $g$ from Proposition \ref{prop:combin}, obeys $$ \|u^{1/a}\|_\infty\lesssim1\,. $$
\end{lemma}

\begin{proof}
We choose $a_1\leq1/4$ and decrease it further if needed.  Lemma~\ref{lem:weakregularity} gives a continuous representative of $u$. Set  $U\coloneqq \max u\geq1$, which is attained at $\theta_0$.

\emph{Step 1: A persisting high peak.}
For $z\in Z(u)$, set $d=1-|z|^2$. Using the disk involution $\varphi_z(w)=(z-w)/(1-\overline z\,w)$, whose Jacobian restricted to $\mathbb S^1$ is $P_z$, gives $\int P_z^{1-a}\,\mathrm d\mu=\int P_z^a\,\mathrm d\mu$. After rotation, expanding  the conjugate factors $(1-|z|e^{i\theta})^{-a}$ and $(1-|z|e^{-i\theta})^{-a}$  into binomial series and integrating their product gives
$$ \int P_z^a\,\mathrm d\mu =d^a\sum_{j\geq0} \left(\frac{\Gamma(a+j)}{j!\Gamma(a)}\right)^2|z|^{2j}
   \leq d^aD_a\,.$$ Here we used that the integral of $e^{i(j-k)\theta}$, $j,k\in \Z$, equals $1$ if $j=k$ and vanishes otherwise. The series has nonnegative coefficients and its value at $|z|=1$ is $\int_{\mathbb S^1} |1-e^{i\theta}|^{-2a} \,\mathrm d\mu=D_a$ by \eqref{eq:inversekernel}. This follows by monotone convergence of the series and dominated convergence of the integral as $2a<1$. Since $\langle u,P_z^{1-a}\rangle=1$, we deduce that
\begin{equation}
     1\leq UD_ad^a\leq Ud^a e^{Ca^2}\,,\label{eq:bubbleboundary}
\end{equation} using \eqref{eq:logD_aDelh} for some $C>0$. Explicit computation gives
$\|P_z\|_\infty\lesssim d^{-1}$, $\|(\log P_z)'\|_\infty\lesssim d^{-1}$, $
\|(\log P_z)''\|_\infty\lesssim d^{-2}
$, and 
$ \|(P_z^a)''\|_\infty \lesssim a d^{-2-a} \lesssim aU^{2/a+1} $ for $0<a\leq 1/4$ by the chain rule, \eqref{eq:bubbleboundary}, and \eqref{eq:logD_aDelh}.
This implies 
\begin{equation}
    \|g''\|_\infty\lesssim aU^{2/a+1}\,.\label{eq:g''}
\end{equation}

Applying
$\Delta_hT_a^{-1}$ to \eqref{eq:combin} yields
$$
(1-\gamma)\Delta_hu+\gamma \Delta_h g
=\|u\|_{1/a}^{2-1/a}(\Delta_hk_a)*u^{1/a-1}\,.
$$
Since $\|u\|_{1/a}^{2-1/a}\leq1$, Lemma~\ref{lem:kernelbounds} and \eqref{eq:g''} imply that
\begin{equation}\label{eq:ineq_chain}
    0\leq U-u(\theta_0+h)\leq -\Delta_h u(\theta_0)\leq C_1a (U^{1/a-1}|h|^{1-2a} +U^{2/a+1}h^2) 
\end{equation} for every $|h|\leq1/10$ and some $C_1>0$. Choose $\eta\in (0,1/10)$ sufficiently small so that $C_1(\eta^{1-2a}+\eta^2)\leq1/2$. 
Since $a\leq1/4$, the preceding bound is at most $aU/2$ whenever $|h|\leq \eta U^{-1/a}$. Then \eqref{eq:ineq_chain} implies
\begin{equation}\label{eq:peakarc} u(\theta_0+h)\geq U(1-a/2)\,.
\end{equation}

\emph{Step 2: Controlling the height.}
Define the width parameter $\ell \coloneqq \eta U^{-1/a}$. The function $u-1$ has mean zero. Using the Fourier expansion without the $n=0$ mode, the Cauchy--Schwarz inequality gives 
\begin{equation*}
 \left|
 \frac1{2\ell}\int_{-\ell}^{\ell}
 (u(\theta_0+t)-1)\,\mathrm dt
 \right|^2
 \lesssim \left(\sum_{n\neq0}\lambda_{|n|}(a)|\widehat {(u-1)}(n)|^2\right)
2a \sum_{n\geq1}n^{2a-1}
 \left(\frac{\sin(n\ell)}{n\ell}\right)^2\,,
\end{equation*} where we used \eqref{eq:symbolbounds}. The first sum is of order $O(a)$ by Lemma~\ref{lem:cap}. Since $t\mapsto t^{2a-1}\min\{1,(t\ell)^{-2}\}$ is
decreasing for $a\leq1/4$, we deduce that
\begin{align*}
2a^2\sum_{n\geq1}n^{2a-1}
 \left(\frac{\sin(n\ell)}{n\ell}\right)^2&\leq 2a^2+2a^2\int_1^\infty t^{2a-1}\min\{1,(t\ell)^{-2}\}\,\mathrm dt\\
&=2a^2+a(\ell^{-2a}-1)+\frac {a^2}{1-a}\ell^{-2a}\lesssim_\eta a(U^2-1) +a^2U^2\,.
\end{align*} Note that $\eta$ is a fixed numerical constant.

The lower bound \eqref{eq:peakarc} gives
$$ 1-\frac 1 U\leq\frac a2+
\frac1U\left|\frac1{2\ell}\int_{-\ell}^\ell (u(\theta_0+t)-1)\,\mathrm dt\right|.
$$ Squaring and using the preceding estimates, $\xi=1-U^{-1}$ satisfies
$$ \xi^2\lesssim a^2+a(1-U^{-2})\leq a^2+2a\xi\,, $$
 so $\xi\leq C_*a$ for some $C_*>0$.
Choose $0<a_1\leq\min\{1/4,(2C_*)^{-1}\}$.
Then $\xi\leq1/2$, so $-\log(1-\xi)\leq 2\xi$ and
$$
\|u^{1/a}\|_\infty
=U^{1/a}
=\exp\!\left(\frac{-\log(1-\xi)}a\right)
\leq e^{2C_*}\,,
$$ which concludes the proof.
\end{proof}

\subsection{Regularity and balancing}
The maximum estimate now makes the nonlinearity in the Euler--Lagrange equation~\eqref{eq:combin} uniformly bounded.  
After subtracting the constant mode, we can apply the inverse of $aT_a$ twice to obtain $\|u-1\|_{H^1}\lesssim a$. This suffices for an expansion and an elementary interpolation argument below. Choosing a conformal map provides the first-moment constraint as required by the limiting inequality.

\begin{proposition}\label{prop:regularity}
There are computable $M>0$ and $a_2\in(0,a_1]$ such that for $0<a\leq a_2$ and $u$ satisfying the assumptions of Lemma \ref{lem:peak}, the critical density $\rho_0$, defined as in \eqref{eq:criticaldensity}, obeys
\begin{equation}\label{eq:Mclass}
 \|\log\rho_0\|_{H^1}\leq M\,.
\end{equation}
After enlarging $M$, the same bound holds for a balanced conformal representative $\rho$ of~$\rho_0$.
\end{proposition}

The critical density $\rho_0$ has unit mass by definition. Balancing preserves this mass but need not preserve \eqref{eq:cap}.

\begin{proof}
    First, we apply  \eqref{eq:bubbleboundary} and Lemma~\ref{lem:peak} to deduce that $$ 1-|z|^2\gtrsim \|u\|_\infty^{-1/a}\gtrsim 1
$$ for every $z\in Z(u)$. Therefore, there is an $R<1$ such that $Z(u)\subset \{|z|<R\} $, so that $T_ag=\int_{Z(u)} P_z^{1-a} \,\mathrm d\tau(z)$, with $g$ as in \eqref{eq:g} and $\tau$ a probability measure, satisfies \begin{equation}
    \label{eq:Tag_bound}
\|T_ag\|_{C^1}\lesssim_R 1\,.\end{equation}

We have $\|u\|_1=1$ by assumption from \eqref{eq:cap}. We now use this bound and the Euler--Lagrange equation \eqref{eq:combin} to bootstrap regularity for $T_a(u-1)$. Indeed, the previous estimates and  \eqref{eq:combin} give $\|T_a(u-1)\|_{\infty}\lesssim 1$. As $u-1$ has mean zero, there is no degree zero mode and \eqref{eq:symbolbounds} implies for $a$ sufficiently small that $n^{4/5}\lesssim a\lambda_n(a)$. This implies a gain of $4/5$ of a derivative in Fourier space when applying the inverse of $aT_a$. In particular, by Morrey embedding, we have \begin{equation*}
    \|u-1\|_\infty\lesssim \|u-1\|_{H^{4/5}}\lesssim a\,.
\end{equation*} Choosing $a$ smaller if necessary, we know that $u\geq 1/2$. By the mean value theorem, we then conclude $\|u^{1/a-1}\|_{H^{1/5}}\lesssim 1+a^{-1}\|u-1\|_{H^{1/5}}\lesssim 1$. Together with \eqref{eq:Tag_bound} and \eqref{eq:combin}, we can upgrade the $L^\infty$-bound on $T_a(u-1)$ to $\|T_a(u-1)\|_{H^{1/5}}\lesssim 1$. By the same gain as before when applying the inverse of $aT_a$, we obtain $$\|u-1\|_{H^1}\lesssim a\,.$$ 
Recall that $1\leq \|u\|_{1/a}^{1/a}\lesssim 1$ by Jensen's inequality and \eqref{eq:upper_crit}. As a consequence of the $H^1$-bound for $(u-1)$, we then have $\|\log(\rho_0)\|_{H^1}\lesssim 1$ for the critical density. This uses $$\log\rho_0=a^{-1}\log u -\log \int_{\mathbb S^1} u^{1/a}\,\mathrm d\mu\qquad \text{and} \qquad (\log\rho_0)'=u' (au)^{-1}\,,$$ and $|\log u|\leq 2|u-1|$. Note that Sobolev embedding also gives $\|\rho_0\|_\infty\lesssim 1$.

Choose a balancing map using \cite[Lemma B.1]{FrankLieb}. After rotation, we can write this map as $\chi_r(\zeta)=(\zeta-r)/(1-r\zeta)$. Balancing and $\int_{\mathbb S^1} \chi_r\cdot e_1\,\mathrm d\mu= -r$ give
$$1=\int_{\mathbb S^1} (1+\chi_r\cdot e_1) \rho_0\,\mathrm d\mu\leq \|\rho_0\|_\infty (1-r)\,.$$ Hence, $r$ stays uniformly below one. For $\phi=\chi_r^{-1}$ and $\rho=J_\phi (\rho_0\circ \phi)$, the chain rule applied to $\log\rho=(\log\rho_0)\circ\phi+\log J_\phi$ gives the same desired bound for $\rho$ after enlarging $M$. 
\end{proof}

\section{The limiting inequality}\label{sec:4}
Proposition \ref{prop:regularity} supplies a balanced probability density $\rho$ satisfying \eqref{eq:Mclass}, obtained by conformally transforming the original critical density. We work with this balanced representative from now on. We next identify a constrained Beckner--Onofri-type inequality as its limit for $a\downarrow0$.

For $u\in H^{1/2-a}$, define the auxiliary functional
\begin{equation}\label{eq:Phi}
 \Phi_a(u)\coloneqq (1-c_{\mathrm{loc}})\langle u,T_au\rangle
+c_{\mathrm{loc}}\left(\int_{\mathbb S^1}u\,\mathrm d\mu\right)^2
 -\|u\|_{1/a}^2\,.
\end{equation} For a probability  density $\rho$, we have $\|\rho^a\|_{1/a}^2=1$. The quotient bound $\mathcal E_a (\rho^a)\leq c_{\mathrm{loc}}$ therefore implies
$\Phi_a(\rho^a)\leq0$, since the constant bubble is an admissible competitor in the supremum.

\subsection{The Beckner--Onofri deficit and its remainder}
The linear term in an expansion of $\Phi_a$ in $a$ turns out to be twice the \textit{constrained Beckner--Onofri deficit}.  For a positive probability density $\rho=e^f$, the latter is given by
\begin{equation}\label{eq:D}D(f)\coloneqq\int_{\mathbb{S}^1} f \,\mathrm d\mu+\frac14 \langle f, (-\Delta_{\mathbb S^1})^{1/2}f\rangle\,.
\end{equation} Here  the operator $(-\Delta_{\mathbb S^1})^{1/2}$ acts via multiplication by $|n|$ on the basis of eigenfunctions $e^{in\theta}$. 
The quantity $D(\log \rho)$ vanishes for the \emph{comparison densities}
\begin{equation}\label{eq:Qdefinition}
 Q_z(\theta)\coloneqq P_z(2\theta)\,,\qquad z\in \mathbb D\,, \end{equation}  where $P_z$ is the Poisson kernel defined in \eqref{eq:bubbles}.  Since $\int_{\mathbb{S}^1}Q_z(\theta)\omega(2\theta)\,\mathrm d\mu=z$, the unique parameter matching the second Fourier moment of a strictly positive balanced probability density $\rho$ is $z=\int_{\mathbb{S}^1}\rho(\theta)\omega(2\theta)\,\mathrm d\mu$, which satisfies $|z|<1$ and $\widehat Q_z(2)=z_1-iz_2$. 

For an integrable density~$\rho$, define its \emph{Toeplitz matrix} by
\begin{equation*}
T_N(\rho)\coloneqq (\widehat\rho(j-k) )_{j,k=0}^{N}\,, \qquad N\in\N\cup\{0\}\,. 
\end{equation*} The next lemma shows how a $3\times 3$ Toeplitz determinant controls both the previous deficit and the distance to its equality family.  

\begin{lemma}\label{lem:szego}
Let $\rho>0$ be a balanced probability density with $\log\rho\in H^1$, and let $Q=Q_z$ have the same second Fourier moment as $\rho$. Then
\begin{equation}\label{eq:entropy}
 \int_{\mathbb{S}^1} (\rho/Q-1-\log(\rho/Q) )\,\mathrm d\mu
 \leq2D(\log\rho)\,.
\end{equation}
In particular, $D(\log\rho)\geq 0$ and equality holds if and only if  $\rho=Q_z$.
\end{lemma}

\begin{proof}
For the sake of simplicity, we denote $f=\log\rho$. Since~$f\in H^1\subset L^\infty$, the density is bounded from above and away from zero; moreover, we also have $\langle f,(-\Delta_{\mathbb S^1})^{1/2}f\rangle <\infty$.
The normalized Toeplitz determinants are increasing in $N$ and tend to the  strong Szeg\H{o} limit \cite[Theorems 1.1 and 2.4]{Simon}. This gives
\begin{equation}\label{eq:toeplitzbound} \log\det T_N(\rho)
 \leq (N+1)\widehat f(0)+\frac12\langle f,(-\Delta_{\mathbb S^1})^{1/2}f\rangle\,.
\end{equation}
 Whereas normalization and balancing give $\widehat \rho(0)=1$  and $\widehat \rho(1)=0$, the matching second Fourier coefficient satisfies $|\widehat\rho(2)|=|z|<1$. Overall, we obtain
\begin{equation*}T_2(\rho)=\begin{pmatrix}1&0&\overline{\widehat\rho(2)}\\0&1&0\\\widehat\rho(2)&0&1\end{pmatrix}
 \qquad  \text{with}\quad\det T_2(\rho)= 1-|z|^2\,.
\end{equation*}
Since the reciprocal of $Q$ only contains the spherical harmonics of degree $0$ and $2$, moment matching gives
$\int_{\mathbb{S}^1}\rho/Q\,\mathrm d\mu=\int_{\mathbb{S}^1}Q/Q\,\mathrm d\mu=1$.
We can also show that $\int_{\mathbb{S}^1}\log Q\,\mathrm d\mu=\log(1-|z|^2)$ as the logarithm of its denominator has mean zero. Using~\eqref{eq:toeplitzbound} with $N=2$, we thus obtain
\begin{align*}
 \int_{\mathbb{S}^1}(\rho/Q-1-\log(\rho/Q))\,\mathrm d\mu&=\log(1-|z|^2)-\widehat f(0)
 \\&
 \leq 2\widehat f(0)+\frac12\langle f,(-\Delta_{\mathbb S^1})^{1/2}f\rangle
 =2D(f)\,.
\end{align*}
Equality $D(f)=0$ forces $\rho/Q=1$ almost everywhere. Conversely, inserting the Fourier series of $\log Q_z$ leads to $D(\log Q_z)=0$.
\end{proof}

\begin{remark}
The coefficient $1/4$ in \eqref{eq:D} is sharp. Indeed, for $Q_r(\theta)=P_r(2\theta)$ with $0<r<1$, we have $$\int_{\mathbb S^1} \log Q_r \,\mathrm d\mu =\log(1-r^2)\qquad \text{and} \qquad \langle \log Q_r, (-\Delta_{\mathbb S^1})^{1/2}\log Q_r\rangle =-4\log(1-r^2)\,.$$ 
Moreover, 
the above equalities also show that a distance estimate must retain the whole family $\{Q_r\}$; the balance assumption is essential, since the ordinary Poisson density satisfies $D(\log P_r)=\frac12\log(1-r^2)<0$.
\end{remark}

We call $w$ \emph{even} if $w(\theta+\pi)=w(\theta)$ and $k$ \emph{odd} if $k(\theta+\pi)=-k(\theta)$. We remark that the comparison densities $Q_z$ are even in the above sense.  For odd perturbations, the mass and second moments vanish automatically, so only the first moment needs to be constrained.
\begin{lemma}\label{lem:oddnormal}
Let $Q=Q_z$, $z\in\mathbb D$, be defined by~\eqref{eq:Qdefinition}. Every odd $k\in H^{1/2}$ satisfying $\int_{\mathbb S^1} Qk\omega\,\mathrm d\mu=0$ obeys
\begin{equation}\label{eq:oddnormalgap}
 \frac12\langle k,(-\Delta_{\mathbb S^1})^{1/2}k\rangle-\int_{\mathbb S^1}Qk^2\,\mathrm d\mu
 \geq\frac12\|k\|_2^2\,.
\end{equation}

\end{lemma}

\begin{proof}
First let $k$ be smooth. For small $|t|$, the curve $$\rho_t=Q(1+tk)$$ consists of densities, which are positive, balanced, and have the same mass and  second moment as $Q$. Note that $\partial_t\log(\rho_t)|_{t=0}=k$, $\partial_t^2\log(\rho_t)|_{t=0}=-k^2$, $\int Qk\,\mathrm d\mu=0$, and $(-\Delta_{\mathbb S^1})^{1/2}\log Q=2(Q-1)$ by comparing Fourier coefficients when expanding the denominator of $Q$. These formulas and direct computation show that
\begin{align*}
    \frac{\mathrm d^2}{\mathrm dt^2} D(\log(\rho_t))\bigg|_{t=0} &=\frac12\langle k,(-\Delta_{\mathbb S^1})^{1/2}k\rangle -\int_{\mathbb S^1} Q k^2\,\mathrm d\mu\,,
\\ \frac{\mathrm d^2}{\mathrm  dt^2}
\int_{\mathbb S^1}\log\frac{Q}{\rho_t}\,\mathrm d\mu\bigg|_{t=0}  &=\|k\|_2^2\,, \end{align*} while first derivatives vanish by parity. Applying \eqref{eq:entropy} and using equality at $t=0$ proves~\eqref{eq:oddnormalgap}.

For general $k$, we can approximate $k$ by odd, smooth Fourier polynomials $k_N$ for $N\in \N$.  We correct the first moments of $k_N$ by linear combinations of
$\cos\theta$ and $\sin\theta$. Since $Q>0$, the matrix $\int_{\mathbb S^1} Q\omega\otimes\omega \,\mathrm d\mu$ is positive definite. Thus, there are $c_N\in \R^2$ such that $\tilde k_N=k_N-c_N\omega$ satisfies $\int_{\mathbb S^1} Q\tilde k_N\omega\,\mathrm d\mu=0$. Note that $c_N\to 0$ as $k_N\to k$ in $H^{1/2}$; we also have $\tilde k_N\to k$ in $H^{1/2}$. Passing to the limit in the inequality completes the proof in the general setting.
\end{proof}

\begin{remark}
At $Q=1$ the odd functions $k=\cos\theta$ and $k=\sin\theta$ make the left side of \eqref{eq:oddnormalgap} equal to $-\|k\|_2^2/2$.  Thus, first-moment constraints are necessary for Lemma~\ref{lem:oddnormal}. Moreover, the admissible function $k=\cos 3\theta$ gives equality in \eqref{eq:oddnormalgap}, so the coefficient $1/2$ is sharp.
\end{remark}

\subsection{Localization near the equality family}
A nonpositive value of $\Phi_a$ forces the limiting deficit to be $O(a)$. The entropy remainder then controls the squared $L^2$-distance to the closest $Q_z$ with matching second moment. The uniform $H^1$-bound upgrades this to uniform proximity. A positive lower bound for $\rho$ keeps $|z|$ away from one. This is the content of the next proposition.

\begin{proposition}\label{prop:localization}
For every $M>0$, there are computable $C_M\geq1$ and $R_M\in(0,1)$ such that the following holds.
Let $\rho$ be a positive balanced probability density with
$\|\log\rho\|_{H^1} \leq M$.
If $0<a\leq 1/4$ and $\Phi_a(\rho^a)\leq0$, then the density $Q_z$ with the same second Fourier moment
as $\rho$ satisfies \begin{equation}\label{eq:localization} |z|\leq R_M \qquad  \text{and}\qquad \|\rho/Q_z-1\|_\infty \leq C_M a^{1/4}\,.
\end{equation}
\end{proposition}

\begin{proof} 
We divide the proof into several steps for the sake of clarity.

\emph{Step 1: Bounding $Q_z$.} 
The assumption $\|\log(\rho)\|_{H^1}\leq M$ implies that $\rho$ and $\rho^{-1}$ are both bounded by some $N_M>1$. Since $\rho$ and $Q_z$ have matching second moments, we deduce that $$|z|=\left|\int(\rho-N_M^{-1})\omega(2\theta)\,\mathrm d\mu\right|
\leq1-N_M^{-1}\,.$$ Hence, we can take $R_M=1-N_M^{-1}$, and as the admissible $z$ are bounded in a compact set, we deduce that $Q_z$ and $Q_z^{-1}$ are bounded uniformly in $C^1$, depending only on $M$.

\emph{Step 2: A bound for the limiting deficit.}
 Rewriting the factors in the product formula~\eqref{eq:gammaprod}
gives $$\frac{a\lambda_n}{n}
=(1-a)\prod_{j=1}^{n-1}\frac{j}{j+1}
\frac{j+1-a}{j+a}\geq (1-a)\prod_{j=1}^{n-1}\left(1-a\left(\frac1j+\frac1{j+1}\right)\right)\,.$$ The inequality $1-\prod_j(1-x_j)\leq\sum_jx_j$ for $x_j\in [0,1]$ and \eqref{eq:symbolbounds} give
\begin{equation}\label{eq:symbolerror}
0\leq n-a\lambda_n
\leq an\left(1+\sum_{j=1}^{n-1}
\left(\frac1j+\frac1{j+1}\right)\right)
\leq2an^2\,.
\end{equation} In the last estimate, we just bound all fractions by $1$. By Taylor's theorem, a short computation using the assumed $H^1$-bound yields $$
\|\rho^a-1-a\log\rho\|_{H^1}\lesssim_M a^2 \qquad \text{and} \qquad  \|\rho^a-1\|_{H^1}\lesssim_M a\,. $$ By~\eqref{eq:symbolerror}, replacing $\lambda_{|n|}$
by $|n|/a$ in Fourier space costs 
\begin{equation*} \left|\sum_{n\neq0} \lambda_{|n|}  |\widehat{\rho^a}(n)|^2 -\frac1a\sum_{n\neq0}|n|\,|\widehat{\rho^a}(n)|^2\right|\leq 2\sum_{n\neq0} n^2  |\widehat{\rho^a}(n)|^2\leq 2\|(\rho^a)'\|_2^2\lesssim_M a^2 \,.
\end{equation*} Hence, we deduce $\left|\sum_{n\neq0} \lambda_{|n|}  |\widehat{\rho^a}(n)|^2-a\langle\log\rho,(-\Delta_{\mathbb S^1})^{1/2}\log {\rho}\rangle\right| \lesssim_M a^2$. Recall the form of $\Phi_a$ in \eqref{eq:Phi}.
Using $\|\rho\|_1=1$ and $1-c_{ \mathrm{loc}}=1/2+ O(a)$, we obtain \begin{align} \notag|\Phi_a(\rho^a)&- 2aD(\log\rho) |\\&=\left|\left(\int_{\mathbb S^1}\rho^a\,\mathrm d\mu\right)^2-1+(1-c_{\mathrm{loc}})\sum_{n\neq0}\lambda_{|n|}|\widehat{\rho^a}(n)|^2 -2aD(\log\rho)\right|\lesssim_Ma^2\,.
\label{eq:PhiatoD}
\end{align}

\emph{Step 3: From entropy to uniform closeness.} Recall from the proof of Lemma \ref{lem:szego} that matching moments give
$\|\rho/Q_z\|_1 = 1 $. The preceding estimates give $\|\rho/Q_z\|_{H^1} \lesssim_M1$
and bound $\rho/Q_z$ above and away from zero. By Lemma \ref{lem:szego}, we deduce that \begin{equation*}
\left\|\frac{\rho}{Q_z}-1\right\|_2^2 \lesssim_M\int_{\mathbb S^1}\left(\frac{\rho}{Q_z}-1 -\log\frac{\rho}{Q_z}\right)\,\mathrm d\mu\leq 2D(\log\rho)\lesssim_M a\,.
\end{equation*}
In the last step, we used \eqref{eq:PhiatoD} and the assumption $\Phi_a(\rho^a)\leq0$. Finally, we bound $\|\rho/Q_z-1\|^2_\infty$ via the one-dimensional Agmon inequality, using $\|\rho/Q_z-1\|_2\lesssim_M\sqrt{a}$ and $\|\rho/Q_z-1\|_{H^1}\lesssim_M1$; this proves \eqref{eq:localization}.
\end{proof}

\section{Comparison by antipodal symmetry}\label{sec:5}
The localization in Proposition~\ref{prop:localization} places the density near an even optimizing density $Q_z$. We first prove nonnegativity for every even function by doubling the angle. Indeed, writing an even function as
$w(\theta)=\widetilde w(2\theta)$ converts its even Fourier
frequencies into all integer frequencies. An eigenvalue
comparison then reduces the desired estimate to the ordinary
sharp Sobolev inequality for $\widetilde w$. We then show that adding the odd part increases the functional near this family. Note that parity makes the first variation in an odd direction at an even function vanish.

\subsection{Nonnegativity for even functions}

Since $P_0=1$, the quotient~\eqref{eq:circleQ} satisfies
\begin{equation*}  
0\geq(\mathcal E_a(u)-c_{\mathrm{loc}} ) (\langle u, T_a u\rangle-\sup_{z\in\mathbb D}|\langle u,P_z^{1-a}\rangle|^2 )\geq\Phi_a(u)
\end{equation*}for $u\notin\mathcal M$ with $\mathcal E_a(u)\leq c_{\mathrm{loc}}$. 
Thus, if we prove that $\Phi_a$ is positive, we obtain a contradiction to the above estimate and thereby exclude global minimizers at or below~$c_{\mathrm{loc}}$.

\begin{lemma}\label{lem:even}
For $0<a<1/2$ and every even $w\in H^{1/2-a}$, we have
\begin{equation}\label{eq:evenpositivity}
 \Phi_a(w)\geq0\,.
\end{equation}
Moreover, equality holds if and only if $w$ is constant.
\end{lemma}
\begin{proof}
We introduce the notation $w(\theta)=\widetilde w(2\theta)$. Thus, $\|w\|_{1/a}=\|\widetilde w\|_{1/a}$, and the means are equal. The definition of Gamma functions allows us to write the recurrence formula $(1-c_{\mathrm{loc}})\lambda_{2}=\lambda_1$ and
$$
\frac{\lambda_{2(n+1)}}{\lambda_{n+1}} \frac{\lambda_{n}}{\lambda_{2n}}
=1+\frac{a(1-a)(1-2a)}
{(2n+a)(2n+1+a)(n+1-a)}>1\,.
$$
Therefore, we obtain $(1-c_{\mathrm{loc}})\lambda_{2n}\geq\lambda_n$,
with strict inequality for $n\geq2$. Comparing Fourier
coefficients and applying the fractional Sobolev inequality \eqref{eq:ineq}, we find
\begin{align*}
(1-c_{\mathrm{loc}})\langle w,T_aw\rangle
+c_{\mathrm{loc}}\left(\int_{\mathbb S^1} w\,\mathrm d\mu\right)^2&=|\widehat{\widetilde w}(0)|^2
+\sum_{n\neq0}(1-c_{\mathrm{loc}})\lambda_{2|n|}|\widehat{\widetilde w}(n)|^2\\&\geq\langle\widetilde w,T_a\widetilde w\rangle\geq\|\widetilde w\|_{1/a}^2
=\|w\|_{1/a}^2\,.
\end{align*}
This proves~\eqref{eq:evenpositivity}.

If equality holds and $w\neq0$, then $\widetilde w$ is a Sobolev optimizer with no Fourier modes of degree at least two. Hence,  we must have $\widetilde w=tP_z^a$. If $z\neq0$, rotation and the binomial expansion of the denominator show that its degree-two coefficient
is nonzero, which is not possible. Thus, $z=0$, so $w$  is constant. Conversely, every constant gives equality.
\end{proof}

\subsection{A spectral bound for odd perturbations}
The density $Q=Q_z$ is even. We extend Lemma \ref{lem:oddnormal} to positive $a$ while allowing an error in the first moments.

\begin{lemma}\label{lem:oddgap}
For every $0<R<1$ there are computable $\varepsilon_R>0$, $\kappa_R>0$, and $C_R>0$ such that, for $Q=Q_z$ with $|z|\leq R$,  $0<a\leq\varepsilon_R$, and every odd $k\in H^{1/2-a}$, we have
\begin{equation}\label{eq:oddgap}
 a(1-c_{\mathrm{loc}}) \langle k,T_ak\rangle-\int_{\mathbb S^1} Qk^2 \,\mathrm d\mu \geq\kappa_R\|k\|_2^2 -C_R\left|\int_{\mathbb S^1} Qk\,\omega\,\mathrm d\mu\right|^2\,.
\end{equation}
\end{lemma}

\begin{proof}
    We first show that there is a computable $K_R>0$ such that
\begin{equation}\label{eq:oddpenalty}
\frac 12\langle k,(-\Delta_{\mathbb S^1})^{1/2}k\rangle -\int_{\mathbb S^1} Qk^2\,\mathrm d\mu
+K_R \left|\int_{\mathbb S^1} Qk\omega\,\mathrm d\mu \right|^2
\geq\frac 1 4\|k\|_2^2\end{equation}
for every odd $k\in H^{1/2}$ and then derive the extension \eqref{eq:oddgap} for odd $k\in H^{1/2-a}$. 
 Let $k_1$ be the $L^2$-projection of $k$ onto
$\operatorname{span}\{Q\cos\theta,Q\sin\theta\}$,
and set $k_0=k-k_1$. Then
$\int_{\mathbb S^1} Qk_0\omega\,\mathrm d\mu=0$. Recall that the Poisson formula gives uniform bounds for $\|Q\|_{C^1}$ and  $\|Q^{-1}\|_\infty$ when $|z|\leq R$ for some $R<1$.
These uniform bounds imply $$\left\| \left(\frac12(-\Delta_{\mathbb S^1})^{1/2 }-Q\right)k_1 \right\|_2\lesssim_R\|k_1\|_2\qquad \text{and}
\qquad \|k_1\|_2\lesssim_R \left|\int_{\mathbb S^1} Qk\omega\,\mathrm d\mu\right|. $$
Applying Lemma \ref{lem:oddnormal} to $k_0$ and
Young's inequality to the cross term gives 
\begin{align*}\frac 12\langle k,(-\Delta_{\mathbb S^1})^{1/2}k\rangle
-\int_{\mathbb S^1} Qk^2\,\mathrm d\mu &\geq\frac 12\|k_0\|_2^2 -C_R(\|k_0\|_2\|k_1\|_2+\|k_1\|_2^2) \\ &\geq\frac 14\|k\|_2^2 -C_R\left|\int_{\mathbb S^1} Qk\omega\,\mathrm d\mu\right|^2,
\end{align*} which proves \eqref{eq:oddpenalty}.

We next treat low and high frequencies separately by splitting $$k=k^{\lo}+k^{\hi}$$ into the Fourier modes with $|n|\leq N$
and with $|n|>N$ for some $N\in\N$ that we choose later. For $0<a\leq1/4$ and $n\geq 1$, estimates \eqref{eq:symbolbounds}
and \eqref{eq:symbolerror} lead to
\begin{equation}
    \label{eq:twobounds}a(1-c_{\mathrm{loc}})\lambda_n\geq\frac14\sqrt n \qquad \text{and} \qquad \left|a(1-c_{\mathrm{loc}})\lambda_n-\frac n2\right| \leq 2 an^2\,. 
\end{equation}
For the low frequencies, the second bound in \eqref{eq:twobounds} and
\eqref{eq:oddpenalty} imply $$ a(1- c_{\mathrm{loc}})\langle k^{\lo},T_ak^{\lo}\rangle
-\int_{\mathbb S^1} Q(k^{\lo})^2\,\mathrm  d\mu +K_R\left|\int_{\mathbb S^1} Qk^{\lo}\omega\,\mathrm d\mu\right|^2 \geq\left( \frac14-2aN^2\right)\|k^{\lo}\|_2^2. $$ For the high frequencies, the first bound in \eqref{eq:twobounds} gives
$$ a (1-c_{\mathrm{loc}})\langle k^{\hi},T_a k^{\hi}\rangle
-\int_{\mathbb S^1} Q(k^{\hi})^2\,\mathrm d\mu \geq\left (\frac14\sqrt N-\|Q\|_\infty\right)\|k^{\hi}\|_2^2\,. $$ For the mixed terms, orthogonality between the modes eliminates
the $T_a$ contribution, while Young's inequality handles
$$
2\left|
\int_{\mathbb S^1} Qk^{\lo} k^{\hi}\,\mathrm d\mu
-K_R\left(\int_{\mathbb S^1}  Qk^{\lo}\omega\,\mathrm d\mu\right)
  \cdot\left(\int_{\mathbb S^1} Qk^{\hi}\omega\,\mathrm d\mu\right) \right| \leq\frac1{16}\|k^{\lo}\|_2^2+
C_R\|k^{\hi}\|_2^2. $$
Combining these bounds and enlarging $C_R$ if necessary, we obtain 
\begin{align*} 
a(1-c_{\mathrm{loc}})\langle k,T_ak\rangle -\int_{\mathbb S^1} Qk^2\,\mathrm d\mu
&+K_R\left|\int_{\mathbb S^1} Qk\omega\,\mathrm d\mu\right|^2\\ 
&\geq \left(\frac3{16}-2aN^2\right) \|k^{\lo}\|_2^2 +\left(\frac14\sqrt N-C_R\right)\|k^{\hi}\|_2^2. 
\end{align*}
We then choose $N$ such that $\sqrt N/4-C_R\geq1/16$ and $\varepsilon_R$ such that $\varepsilon_R\leq1/4$ and $\varepsilon_R N^2\leq1/16$.
By construction, both coefficients are at least $1/16$, which proves
\eqref{eq:oddgap}. Since we used \eqref{eq:oddpenalty} only for the smooth low-frequency part, the argument is easy to extend to $H^{1/2-a}$, which completes the proof.
\end{proof}

\subsection{Removing the odd part}
The first moments in~\eqref{eq:oddgap} are linear in the perturbation, whereas balancing provides constraints for $u^{1/a}$. Thus, along the segment joining $u$ to its antipodal translate, the fundamental theorem of calculus converts balancing into an exact weighted first-moment identity with a weight close to $Q$, which makes the penalty in~\eqref{eq:oddgap} small.

\begin{proposition}\label{prop:finite}
For every $0<R<1$ there are computable $a_R>0$, $\delta_R>0$, and $b_R>0$ with the following property. Let $Q=Q_z$ with $|z|\leq R$, $0<a\leq a_R$, and $\rho$ be a positive balanced probability density such that $u=\rho^a\in H^{1/2-a}$ and $\|\rho/Q-1\|_\infty\leq\delta_R$. If 
\begin{equation*}w(\theta)=\frac12 (u(\theta)+u(\theta+\pi) ) \quad\text{and}\quad v(\theta)=\frac12 (u(\theta)-u(\theta+\pi) )\,.
\end{equation*}
denote the (antipodally) even and odd parts of $u$, respectively, then
\begin{equation}\label{eq:finitecomparison}
 \Phi_a(u)\geq\Phi_a(w)+\frac{b_R}{a}\|v\|_2^2\geq0\,.\end{equation}
Moreover, $\Phi_a(u)=0$ if and only if $\rho=1$.
\end{proposition}

\begin{proof}
We divide the proof into several steps for the sake of clarity.

\emph{Step 1: Controlling the first moments.} Set $\delta=\|\rho/Q-1\|_\infty$, and
$u_t=w+tv$ for  $0\leq t\leq1$. Assume that $a\leq\min\{\varepsilon_R,1/4\}$ and $\delta\leq1/2$, where $\varepsilon_R$ is obtained in Lemma~\ref{lem:oddgap}. Since $Q$ is even, we deduce by definition of $\delta$ and $u_t$ that $$
    (1-\delta)Q \leq u_t^{1/a}\leq(1+\delta)Q $$ for all $t\in [0,1]$.  The uniform upper and lower bounds for $Q$ consequently give \begin{equation}\label{eq:two_approx}
        |u_t^{1/a-1}-Q|\lesssim_R a+\delta\qquad 
\text{and} \qquad(1-a)\left(\int u_t^{1/a}\,\mathrm d \mu\right)^{2a-1} u_t^{1/a-2}\leq Q+C_R(a+\delta) 
    \end{equation} for some $C_R>0$ uniformly in $t$ and $\theta$. Balancing of $u^{1/a}$ and evenness of $w^{1/a}$ imply $$ 0=a \int(u^{1/a}-w^{1/a})\omega\,\mathrm d\mu =\int_0^1\int u_t^{1/a-1}v\omega\,\mathrm d\mu\,\mathrm dt\,.
$$ Therefore, the first approximation in \eqref{eq:two_approx}  gives \begin{equation}\label{eq:approxoddmoment} \left|\int Qv\omega\,\mathrm d\mu\right|
\lesssim_R(a+\delta)\|v\|_2\,. \end{equation}

\emph{Step 2: Controlling the nonlinear norm.}
As ${1/a}>2$, differentiation and the second approximation in \eqref{eq:two_approx} yield a $C_R>0$ such that
\begin{align*}
\frac a2\frac{\mathrm d^2}{\mathrm dt^2}\|u_t\|_{1/a}^2 &\leq (1-a)\left(\int u_t^{1/a}\,\mathrm d\mu \right)^{2a-1} \int u_t^{{1/a}-2}v^2\, 
     \mathrm  d\mu\\ &\leq\int Qv^2\,\mathrm d\mu
+C_R(a+\delta)\|v\|_2^2\,. \end{align*}     Since the first derivative vanishes at  $t=0$ by parity, Taylor's formula leads to
\begin{equation}\label{eq:normparitycost} a(\|u\|_{1/a}^2-\|w\|_{1/a}^2) \leq\int Qv^2\,\mathrm d\mu
   +C_R(a+\delta)\|v\|_2^2\,.
\end{equation}

\emph{Step 3: Comparing the energies.}
Due to parity we find
$$ \Phi_a(u)-\Phi_a(w) =(1-c_{\mathrm{loc}})\langle v,T_av\rangle -\|u\|_{1/a}^2+\|w\|_{1/a}^2\,. $$ Combining the preceding estimates \eqref{eq:approxoddmoment}  and \eqref{eq:normparitycost} with
Lemma \ref{lem:oddgap}, and using
$(a+\delta)^2\leq a+\delta$, gives $$ a(\Phi_a(u)-\Phi_a(w)) \geq(\kappa_R-C_R(a+\delta))\|v\|_2^2\,.
$$ Choose computable $a_R,\delta_R>0$ so that the coefficient on the right side is at least $b_R=\kappa_R/2$, which, together with Lemma \ref{lem:even}, implies \eqref{eq:finitecomparison}. Equality forces $v=0$ and $w$ constant, so normalization and positivity then give $\rho=1$. The converse is immediate, so we omit it here.
\end{proof}

\section{Completion of the proof}\label{sec:6}

Finally, we prove the two parts of Theorem \ref{thm}.

\subsection{Nonattainment}
We choose the constants before applying the comparison. First, take $a_{2}$ and $M$ from Proposition \ref{prop:regularity}, then $C_M$ and $R_M$ from Proposition~\ref{prop:localization}, and finally $a_R$ and $\delta_R$ with  $R=R_M$ from Proposition~\ref{prop:finite}. We define
\begin{equation*} a_*\coloneqq\frac12\min\left\{\frac1{4},a_{2},a_R,(\delta_R/C_M)^4\right\}\qquad \text{and}\qquad
 s_0\coloneqq\frac12-a_*\,. \end{equation*}
Every entry is positive and has a computable positive lower bound, so $0<s_0<1/2$.

Suppose that $0<a<a_*$ and that a global minimizer outside $\mathcal M$ exists. By the local upper bound~\eqref{eq:localconstant}, its quotient is at most $c_{\mathrm{loc}}$. Proposition \ref{prop:combin} gives a nonnegative minimizer satisfying \eqref{eq:cap} and the Euler--Lagrange equation \eqref{eq:combin}. Lemmas \ref{lem:cap} and~\ref{lem:peak} and Proposition~\ref{prop:regularity} produce a balanced probability density~$\rho$ satisfying \eqref{eq:Mclass}.

The function~$\rho^a=J_\phi^a(u\circ\phi)/\|u\|_{1/a}$ is obtained from the normalized minimizer by a conformal change and multiplication by a constant, so its quotient is unchanged and $\|\rho\|_1=\|\rho^a\|_{1/a}=1$.  Consequently,
\begin{equation}\label{eq:lowdensity}
 \frac{1+a}{2-a}\langle\rho^a,T_a\rho^a\rangle +\frac{1-2a}{2-a}\sup_{z\in\mathbb D}|\langle\rho^a,P_z^{1-a}\rangle|^2\leq1\,.
\end{equation}
As $z=0$ is admissible, we can bound the supremum by $\|\rho\|_a^{2a}$ from below, which gives $\Phi_a(\rho^a)\leq0\,.$

Proposition~\ref{prop:localization} selects a matching~$Q_z$ with~$|z|\leq R_M$ and
$$\|\rho/Q_z-1\|_\infty\leq C_Ma^{1/4}\leq\delta_R\,.$$
All hypotheses of Proposition~\ref{prop:finite} now hold. It gives $\Phi_a(\rho^a)>0$ unless $\rho=1$. The first alternative contradicts~\eqref{eq:lowdensity}, while the second one implies $\rho^a=1$ is a Sobolev optimizer, so the original minimizer lies in $\mathcal M$ after multiplication by a constant and a conformal change. Thus, both alternatives are impossible.

\subsection{The best stability constant}
We use the following implication from K\"onig's compactness proof
\cite[Introduction and Section 4]{KoenigExistence}. At fixed $s$, strict separation from both the local and the two-profile thresholds implies attainment.  In dimension one, the required inequalities are
$$c_{\mathrm{BE}}(s)<c_{\mathrm{loc}}=\frac{4s}{3+2s} \qquad\text{and}\quad c_{\mathrm{BE}} (s)<2-2^{1-2s}\,, $$  where the second inequality follows from the first, since the elementary bound $2^t\leq 1+t$ for $t\in[0,1]$, which follows from convexity of $t\mapsto 2^t$, gives $2-2^{1-2s}\geq 2s > 4s/(3+2s)$ for $0<s<1/2$. Although \cite{KoenigExistence} only considers $d\geq 2$, the paper  explicitly identifies strict separation from the local threshold as the only use of the dimension restriction.

If~$c_{\mathrm{BE}}(s)<c_{\mathrm{loc}}$ in the interval of the theorem, compactness would therefore produce a minimizer, contradicting the preceding subsection. This implies $c_{\mathrm{BE}}(s)\geq c_{\mathrm{loc}}$. The local test sequence in \cite{CFW} gives the reverse inequality, so $c_{\mathrm{BE}}(s)=c_{\mathrm{loc}}$.  This proves the second part of Theorem~\ref{thm}.

\renewcommand{\baselinestretch}{1}\selectfont
\providecommand{\etalchar}[1]{$^{#1}$}

\end{document}